\documentstyle[twoside]{article}
\title{\bf Asymptotics of a Gauss hypergeometric function with large parameters, IV: A uniform expansion}
\author{\sc R. B. Paris\footnote{E-mail address:\ \ {\tt r.paris@abertay.ac.uk}}\\
\\
{\em Division of Computing and Mathematics,}\\
{\em Abertay University, Dundee DD1 1HG, UK}\\
}

\begin{document}
\newcommand{\bee}{\begin{equation}}
\newcommand{\ee}{\end{equation}}
\def\f#1#2{\mbox{${\textstyle \frac{#1}{#2}}$}}
\def\dfrac#1#2{\displaystyle{\frac{#1}{#2}}}
\newcommand{\fr}{\frac{1}{2}}
\newcommand{\fs}{\f{1}{2}}
\newcommand{\g}{\Gamma}
\newcommand{\br}{\biggr}
\newcommand{\bl}{\biggl}
\newcommand{\ra}{\rightarrow}
\renewcommand{\topfraction}{0.9}
\renewcommand{\bottomfraction}{0.9}
\renewcommand{\textfraction}{0.05}
\newcommand{\mcol}{\multicolumn}
\date{}
\maketitle
\pagestyle{myheadings}
\markboth{\hfill {\it R.B. Paris} \hfill}
{\hfill {\it Hypergeometric function with large parameters} \hfill}
\begin{abstract} 
We consider the uniform asymptotic expansion for the Gauss hypergeometric function
\[F(a+\epsilon\lambda,m;c+\lambda;x),\qquad \lambda\to+\infty\]
for $x<1$ and positive integer $m$ when the parameter $\epsilon>1$ and the constants $a$ and $c$ are supposed finite. When $m=1$, we employ the standard procedure of the method of steepest descents modified to deal with the situation when a saddle point is near a simple pole. It is shown that it is possible to give a closed-form expression for the coefficients in the resulting uniform expansion. The expansion when $m\geq 2$ is obtained by means of a recurrence relation. Numerical results illustrating the accuracy of the resulting expansion are given. 
\vspace{0.4cm}

\vspace{0.3cm}

\end{abstract}

\vspace{0.2cm}

\noindent $\,$\hrulefill $\,$

\vspace{0.2cm}

\begin{center}
{\bf 1. \  Introduction}
\end{center}
\setcounter{section}{1}
\setcounter{equation}{0}
\renewcommand{\theequation}{\arabic{section}.\arabic{equation}}
The Gauss hypergeometric function is defined by
\[F\left(\!\!\begin{array}{c}a, b\\c\end{array}\!;z\!\right)=\sum_{n=0}^\infty\frac{(a)_n (b)_n}{(c)_n n!}\,z^n \qquad(|z|<1)\]
and elsewhere by analytic continuation, where $(a)_n=\g(a+n)/\g(a)=a(a+1)\ldots (a+n-1)$ is the Pochhammer symbol or rising factorial. The asymptotic expansion of this function
for large values of the parameters $a$, $b$ and $c$ and fixed complex $z$ was first considered by Watson \cite{W} in 1918  and recently by the author in \cite{P1, P2}; see also \cite{Cvit} for the case of two large parameters.  

In \cite[Section 3]{P1}, the expansion of the function 
\bee\label{e11}
F_b(\lambda;z):=F\left(\begin{array}{c}a+\epsilon\lambda, b\\c+\lambda\end{array}\!;z\right)
\ee
with finite $a$, $b$ and $c$ and complex $z$ was considered for $\lambda\ra+\infty$, $\epsilon>0$ using the method of steepest descents applied to suitable integral representations. 
It should first be noted that when $\epsilon<1$ the function $F_b(\lambda;z)$ converges at $z=1$, since the convergence condition $\Re (c-a-b)+(1-\epsilon)\lambda$ will be positive for sufficiently large $\lambda$. When $\epsilon>1$, however, this condition will be broken as $\lambda\to+\infty$ and  $F_b(\lambda;z)$ will not converge at $z=1$. From \cite[p.~267]{Cop}, the dominant behaviour of $F_b(\lambda;z)$ in this limit is
\[\lim_{z\to 1-}F_b(\lambda;z)=\frac{\g(c+\lambda) \g(a+b-c+(\epsilon-1)\lambda)}{\g(a+\epsilon\lambda) \g(b)}\,(1-z)^{c-a-b-(\epsilon-1)\lambda}.\]

The analysis of the asymptotic expansion of $F_b(\lambda;z)$ required different representations for $\epsilon<1$ and $\epsilon>1$.
The resulting expansion when $\epsilon<1$ is given in \cite[(3.6)]{P1} and has the leading behaviour $F_b(\lambda;z)\sim(1-\epsilon z)^{-b}$; this case is also considered in \cite{LP}.
The asymptotic expansion when $\epsilon>1$ is given by \cite[(3.13)]{P1}
\[F_b(\lambda;z)\sim \frac{\g(c+\lambda)\g(1+a-c+(\epsilon-1)\lambda)}{\g(a+\epsilon\lambda)}\hspace{5cm}\]
\bee\label{e12}
\hspace{3cm}\times \frac{\epsilon^{a-\frac{1}{2}+\epsilon\lambda}(\epsilon-1)^{c-a-\frac{1}{2}+(1-\epsilon)\lambda}}{\sqrt{2\pi}\,(1-\epsilon z)^b} \sum_{k=0}^\infty \frac{c_k \g(k+\fs)}{\lambda^{k+\frac{1}{2}} \g(\fs)}
\ee
as $\lambda\to\infty$, where $c_0=1$ and the coefficients $c_1$ and $c_2$ are explicitly stated in \cite[(3.4), (3.5)]{P1}. Application of Stirling's formula
\[\g(a+x)\sim \sqrt{2\pi}\,x^{x+a-\frac{1}{2}} e^{-x}\qquad (x\to+\infty)\]
shows that the factor involving gamma functions appearing in (\ref{e12}) has the behaviour
\bee\label{e13}
\frac{\g(c+\lambda)\g(1+a-c+(\epsilon-1)\lambda)}{\g(a+\epsilon\lambda)}\sim (2\pi\lambda)^\frac{1}{2} \epsilon^{\frac{1}{2}-a-\epsilon\lambda} (\epsilon-1)^{\frac{1}{2}+a-c+(\epsilon-1)\lambda}\qquad(\lambda\to+\infty),
\ee
so that the leading large-$\lambda$ behaviour of $F_b(\lambda;z)$ again reduces to $(1-\epsilon z)^{-b}$; see also \cite{Cvit}.

The above expansion is of Poincar\'e type and is subject to an inconvenient restriction (when $\epsilon>1$) resulting from the requirement that the singularity of the integrand at $1/z$ should lie outside the closed-loop contour of integration. In the case of real $z$ ($=x$) this restriction corresponds to $\epsilon x<1$.
Thus the result presented in (\ref{e12}) leaves unanswered the expansion of the above hypergeometric function for $1/\epsilon\leq x<1$ when $\epsilon>1$ and also that in the neighbourhood of $\epsilon x=1$.

It is the aim in the present paper to derive the expansion of $F_b(\lambda;x)$ that holds uniformly in the neighbourhood of $\epsilon x=1$. This is obtained in the case $b=1$ by the standard procedure of the method of steepest descents modified to deal with the situation when a saddle point is near a simple pole; see, for example, \cite[Section 1.5.1]{PBook}, \cite[p.~356]{Wong} or \cite{J} for a description of the method. We show that it is possible to give a closed-form expression for the coefficients in this uniform expansion. We also give the explicit expression of the first three of these coefficients at the coalescence of the saddle and pole. These asymptotic results are then extended to $b=m$, where $m$ is a positive integer, by means of a recurrence relation.  Numerical results are presented to demonstrate the accuracy of the expansions obtained.
\vspace{0.6cm}

\begin{center}
{\bf 2. \ A uniform expansion of $F_1(\lambda;x)$ for $\lambda\ra+\infty$}
\end{center}
\setcounter{section}{2}
\setcounter{equation}{0}
\renewcommand{\theequation}{\arabic{section}.\arabic{equation}}
We consider the expansion of the hypergeometric function $F_b(\lambda;x)$ defined in (\ref{e11}) as $\lambda\to+\infty$ when $0<x<1$, $\epsilon>1$ and the parameter $b=1$. From \cite[(15.6.2)]{DLMF} we have the integral representation
\bee\label{e21}
F_1(\lambda;x)=\frac{G(\lambda)}{2\pi i} \int_0^{(1+)} \frac{e^{-\lambda\psi(t)}f(t) }{1-xt}\,dt,
\ee
where the integration path is a loop that starts at $t=0$, encircles the point $t=1$ in the positive sense (excluding the point $t=1/x$) and returns to $t=0$, and
\[\psi(t)=(\epsilon-1) \log (t-1)-\epsilon \log\,t,\qquad f(t)=t^{a-1}(t-1)^{c-a-1},\]
with
\[G(\lambda)=\frac{\g(c+\lambda) \g(1+a-c+(\epsilon-1)\lambda)}{\g(a+\epsilon\lambda)}\sim 
\frac{(2\pi\lambda)^\frac{1}{2}\epsilon^{\frac{1}{2}-a}}{(\epsilon-1)^{c-a-\frac{1}{2}}}\,e^{\lambda\psi(\epsilon)} \quad (\lambda\to+\infty)\]
by (\ref{e13}).
The $t$-plane has a branch cut along $(-\infty,1]$ and the integrand has a simple pole at $t=1/x$. The exponential factor has a saddle point where $\psi'(t)=(\epsilon-t)/(t(t-1))=0$; that is, at the point $t_s=\epsilon$, where $\psi''(t_s)=-\epsilon^{-1}(\epsilon-1)^{-1}<0$. The path of steepest descent through the saddle has directions $\pm\fs\pi$ at $t_s$ and forms a closed loop surrounding $t=1$ with endpoints at $t=0$.

To account for the proximity of the saddle point $t_s$ to the pole, we make the standard substitution
\bee\label{e21a}
u^2=\psi(t)-\psi(t_s),
\ee
so that
\bee\label{e22}
F_1(\lambda;x)=\frac{G(\lambda) e^{-\lambda\psi(t_s)}}{2\pi i (-x) }\int_{-\infty}^\infty e^{-\lambda u^2} \frac{f(t)}{t-\alpha}\,\frac{dt}{du}\,du,\qquad \alpha:=1/x.
\ee
The saddle point corresponds to $u=0$ and the pole is now situated at $u=u_\alpha$, where
\bee\label{e22a}
u_\alpha^2=\psi(\alpha)-\psi(t_s).
\ee

The steepest descent path through the saddle $t_s$ maps to the real $u$-axis and
the domain enclosed by this path in the $t$-plane maps into the upper half of the $u$-plane. Thus, we have 
\bee\label{e23e}
u_\alpha=\mp ip\quad\left\{\begin{array}{l}\alpha>\epsilon\ (\epsilon x<1)\\\alpha<\epsilon\ (\epsilon x>1)\end{array}\right.\!\!,\qquad p=(\psi(t_s)-\psi(\alpha))^{1/2}>0.
\ee
In the neighbourhood of the saddle point
\bee\label{e23}
\frac{f(t)}{t-\alpha}\,\frac{dt}{du}=\frac{d_{-1}}{u-u_\alpha}+h(u),
\ee
where $h(u)$ is regular at $u=u_\alpha$ and $u=0$. The coefficient $d_{-1}$ is given by
\bee\label{e23a}
d_{-1}=\mathop{\lim}_{\scriptstyle t\to \alpha \atop \scriptstyle u\to u_\alpha} \bl(\frac{u-u_\alpha}{t-\alpha}\br) f(t) t'(u)=f(\alpha).
\ee

\vspace{0.3cm}

\noindent{\bf 2.1\ \ The expansion of $F_1(\lambda;x)$}.\ \ \ 
When $\alpha>\epsilon$ the pole at $t=\alpha$ is situated outside the loop (corresponding to the steepest descent path) in the integral (\ref{e21}) and $u_\alpha=-ip$. Then
\[F_1(\lambda;x)=\frac{G(\lambda)e^{-\lambda\psi(t_s)}}{2\pi i(-x)}
\bl\{d_{-1}\int_{-\infty}^\infty \frac{e^{-\lambda u^2}}{u+ip}\,du+\int_{-\infty}^\infty e^{-\lambda u^2} h(u)\,du\br\},\] 
where
\bee\label{int}
\int_{-\infty}^\infty \frac{e^{-\lambda u^2}}{u\pm ip}\,du=\mp\pi ie^{\lambda p^2} \mbox{erfc}\,(\lambda^\frac{1}{2}p)
\ee
with erfc denoting the complementary error function. Introduction of the expansion
\bee\label{e24}
h(u)\sim \sum_{k=0}^\infty {\hat d}_k u^k,
\ee
followed by application of Watson's lemma, then yields the expansion as $\lambda\to+\infty$ when $\alpha>\epsilon$
\bee\label{e25}
F_1(\lambda;x)\sim\frac{G(\lambda)}{2x}
\bl\{e^{-\lambda \psi(\alpha)}f(\alpha)\,\mbox{erfc}\,(\lambda^\frac{1}{2}p)+\frac{e^{-\lambda\psi(t_s)}}{\pi}\sum_{k\geq 0} d_{2k}\frac{ \g(k+\fs)}{\lambda^{k+\frac{1}{2}}}\br\},
\ee
where $d_{2k}=i{\hat d}_{2k}$ and from (\ref{e22a}) we have used $\psi(t_s)-p^2=\psi(\alpha)$.

When $\alpha<\epsilon$ we have $u_\alpha=+ip$. Then, employing the result in (\ref{int}), we find the contribution to the integral (\ref{e22}) given by
\[\frac{G(\lambda)e^{-\lambda\psi(\epsilon)}}{2 x}\bl\{
-f(\alpha)\,e^{\lambda p^2}\mbox{erfc}\,(\lambda^\frac{1}{2}p)
+\frac{1}{\pi}\sum_{k\geq 0} d_{2k}\frac{\g(k+\fs)}{\lambda^{k+\frac{1}{2}}}\br\}.\]
In expanding the integration path of the integral (\ref{e21}) to pass through the saddle point $t_s=\epsilon$ it is now necessary pass over the pole at $t=\alpha$, thereby receiving a residue contribution given by $-f(\alpha) e^{-\lambda\psi(\alpha)}$.

Then, upon use of the result $\mbox{erfc}\,(-x)=2-\mbox{erfc}\,(x)$, we obtain the following theorem:
\newtheorem{theorem}{Theorem}
\begin{theorem}$\!\!\!.$\  The
expansion of $F_1(\lambda;x)$ as $\lambda\to+\infty$ is  given by
\bee\label{e27b}
F_1(\lambda;x)\sim\frac{G(\lambda)}{2 x}
\bl\{e^{-\lambda\psi(\alpha)}f(\alpha)\,\mbox{erfc}\,(\pm\lambda^\frac{1}{2}p)+\frac{e^{-\lambda\psi(\epsilon)}}{\pi}
\sum_{k\geq 0} d_{2k}\frac{\g(k+\fs)}{\lambda^{k+\frac{1}{2}}}\br\},
\ee
where the upper sign applies when $\alpha>\epsilon$ $(\epsilon x<1)$ and the lower sign when $\alpha<\epsilon$ $(\epsilon x>1)$. The quantity $p$ is defined in (\ref{e23e}) and vanishes at coalescence. 
\end{theorem}

An alternative form of this expansion is given in the appendix. 
\vspace{0.3cm}

\noindent{\bf 2.2\ \ The coefficients $d_{2k}$}. \ \ \  
The coefficients $d_{2k}$ can be obtained by differentiation of (\ref{e23}) and (\ref{e24}) with respect to $u$ with the derivatives evaluated at $u=0$, $t=\epsilon$. Thus we find
\[d_{2k}=\frac{i}{(2k)!} \bl(\frac{d}{du}\br)^{2k}\bl\{\frac{f(t)t'(u)}{t-\alpha}-\frac{d_{-1}}{u-u_\alpha}\br\}\br|_{u=0,\, t=\epsilon}
\equiv c_{2k}-b_{2k}.\]

The coefficients $c_{2k}$ are given by
\[c_{2k}=\frac{i}{(2k)!} \bl(\frac{d}{du}\br)^{2k}\bl(\frac{f(t)t'(u)}{t-\alpha}\br)\br|_{u=0,\, t=\epsilon},\]
with
\bee\label{e26a}
c_0=\frac{f(\epsilon) t'(0)}{i(\alpha-\epsilon)}=\frac{f(\epsilon)}{\kappa (\alpha-\epsilon)},\qquad \kappa:=(\fs|\psi''(\epsilon)|)^{1/2}
\ee
since by differentiation of (\ref{e21a}) $t'(0)=(2/\psi''(\epsilon))^{1/2}$, where we take\footnote{Since $\psi''(\epsilon)<0$ it is necessary to specify the value of $(\psi''(\epsilon))^{1/2}$. This follows from the fact that in the integral (\ref{e21}) we require that the quantity $(t-t_s)^2 \psi''(\epsilon)/2$ appearing in the exponential be positive on the steepest descent path. Since $\arg (t-t_s)=\fs\pi$ on leaving the saddle then $\arg\,\psi''(\epsilon)=-\fs\pi$.}
$\arg\,\psi''(\epsilon)=-\fs\pi$.

A general expression for the normalised coefficients $C_{2k}=c_{2k}/c_0$ is given by the Wojdylo formula \cite{Woj}
\bee\label{e25w}
C_{2k}=\frac{1}{{\hat\alpha}_0^{k}}\sum_{s=0}^{2k} \frac{{\hat\beta}_{2k-s}}{{\hat\beta}_0} \sum_{j=0}^s \frac{(-)^j (s+\fs)_j}{j!\ {\hat\alpha}_0^j}\,{\cal B}_{kj}\,;
\ee
see also \cite[p.~25]{T}. Here ${\cal B}_{kj}\equiv {\cal B}_{kj}({\hat\alpha}_1, {\hat\alpha}_2, \ldots , {\hat\alpha}_{k-j+1})$ are the partial ordinary Bell polynomials generated by the recursion\footnote{For example, this generates the values ${\cal B}_{41}={\hat\alpha}_4$, ${\cal B}_{42}={\hat\alpha}_2^2+2{\hat\alpha}_1{\hat\alpha}_3$, ${\cal B}_{43}=3{\hat\alpha}_1^2{\hat\alpha}_2$ and ${\cal B}_{44}={\hat\alpha}_1^4$.}
\[{\cal B}_{kj}=\sum_{r=1}^{k-j+1} {\hat\alpha}_r {\cal B}_{k-r,j-1} ,\qquad {\cal B}_{k0}=\delta_{k0},\]
where $\delta_{mn}$ is the Kronecker symbol, and the coefficients ${\hat\alpha}_r$ and ${\hat\beta}_r$ appear in the expansions
\[
\psi(t)-\psi(t_s)=\sum_{r=0}^\infty {\hat\alpha}_r (t-t_s)^{r+2},\qquad {\hat f}(t):=\frac{f(t)}{t-\alpha}=\sum_{r=0}^\infty {\hat\beta}_r(t-t_s)^r
\]
valid in a neighbourhood of the saddle $t=t_s$, where ${\hat\alpha}_r=\psi^{(r+2)}(\epsilon)/(r+2)!$.

Explicit representations for the first three normalised coefficients $C_{2k}$ are
\[C_0=1,\qquad C_2=\frac{1}{\psi''}\{{\hat F}_2-\Psi_3{\hat F}_1+\f{5}{12}\Psi_3^2-\f{1}{4}\Psi_4\},\]
\[C_4=\frac{1}{(\psi'')^2}\{\f{1}{6}{\hat F}_4-\f{5}{9}\Psi_3{\hat F}_3+\f{5}{12}(\f{7}{3}\Psi_3^2-\Psi_4){\hat F}_2-\f{35}{36}(\Psi_3^3-\Psi_3\Psi_4+\f{6}{35}\Psi_5){\hat F}_1 \]
\bee\label{e25c}
+\f{35}{36}(\f{11}{24}\Psi_3^4-\f{3}{4}(\Psi_3^2-\f{1}{6}\Psi_4)\Psi_4+\f{1}{5}\Psi_3\Psi_5-\f{1}{35}\Psi_6)\},
\ee
where, for brevity, we have defined
\bee\label{e25d}
\Psi_k:=\frac{\psi^{(k)}(\epsilon)}{\psi''(\epsilon)}\ \ \ (k\geq 3),\qquad {\hat F}_k:=\frac{{\hat f}^{(k)}(\epsilon)}{{\hat f}(\epsilon)}\ \ \ (k\geq 1);
\ee
see, for example, \cite[p.~119]{D}, \cite[p.~127]{O} or \cite[p.~13]{PBook}. 

The coefficients $b_{2k}$ are given by
\[b_{2k}=\frac{i d_{-1}}{(2k)!} \bl(\frac{d}{du}\br)^{2k} \frac{1}{u-u_\alpha} \br|_{u=0}=\frac{i d_{-1}}{(-u_\alpha)^{2k+1}}\]
\bee\label{e26}
=\pm\frac{(-)^kf(\alpha)}{p^{2k+1}}.
\ee
Then we obtain the closed-form representation for the coefficients $d_{2k}$ appearing in the expansion (\ref{e25}) given by
\bee\label{e27a}
d_{2k}=\frac{f(\epsilon) C_{2k}}{\kappa (\alpha-\epsilon)}\mp\frac{(-)^kf(\alpha)}{p^{2k+1}}\qquad(k=0, 1, 2, \ldots),
\ee
where $\kappa$ and $C_{2k}$ are specified in (\ref{e26a}) and (\ref{e25w}). The upper or lower signs in (\ref{e26}) and (\ref{e27a}) are chosen according as $\epsilon x<1$ or $\epsilon x>1$, respectively. However, these coefficients present a removable singularity when the saddle and pole coincide since both $\alpha-\epsilon$ and $p$ vanish in this limit.
This case is considered in the next section.

\vspace{0.6cm}

\begin{center}
{\bf 3. \  The expansion at coalescence $\alpha=\epsilon$}
\end{center}
\setcounter{section}{3}
\setcounter{equation}{0}
\renewcommand{\theequation}{\arabic{section}.\arabic{equation}}
The expansions in (\ref{e27b}) are suitable when the pole at $t=\alpha$ is not too close to the saddle at $t=\epsilon$. As $\alpha\to\epsilon$ (that is, as $\epsilon x\to 1$) the coefficients $d_{2k}$ present a removable singularity at $t=\epsilon$, $u=0$.

Let $\delta=\alpha-\epsilon$ provide a measure of proximity to coalescence. Then as $\delta\to 0$ the coefficients $c_{2k}$ involve terms of $O(\delta^{-r})$ for integer $r$ in the range $1\leq r\leq 2k+1$. These singular terms cancel with the corresponding terms\footnote{This has been explicitly verified for the coefficients $d_{2k}$ with $k\leq 2$.} present in the coefficients $b_{2k}$ to leave terms of order $O(\delta^r)$, $r=0, 1, 2, \ldots\ $. To see this we make use of the expansions
\[f(\alpha)=f(\epsilon)\{1+F_1\delta+\fs F_2\delta^2+\cdots\},\qquad F_k:=\frac{f^{(k)}(\epsilon)}{f(\epsilon)}\]
and, from (\ref{e22a}),
\[p=\kappa\delta\bl\{1+\frac{2\Psi_3\delta}{3!}+\frac{2\Psi_4\delta^2}{4!}+\cdots\br\}^{1/2}, \]
where $\kappa$ and $\Psi_k$ are defined in (\ref{e26a}) and (\ref{e25d}). Then as $\delta\to 0$
\[b_{2k}=\frac{(-)^k f(\alpha)}{p^{2k+1}}=\frac{(-)^k f(\epsilon)}{(\kappa\delta)^{2k+1}}\,\frac{\{1+F_1\delta+F_2\delta^2/2+\cdots\}}{\{1+2\Psi_3\delta/3!+2\Psi_4\delta^2/4!+\cdots \}^{k+\frac{1}{2}}}\]
\[=\frac{(-)^kf(\epsilon)}{(\kappa\delta)^{2k+1}}\bl\{1+\sum_{s=1}^\infty D_s(k)\delta^s\br\}.\]
The values of the coefficients $d_{2k}$ at coalescence therefore involve the $O(\delta^0)$ term in the above expansion to yield
\bee\label{e28b}
d_{2k}=\frac{(-)^{k-1} f(\epsilon)}{\kappa^{2k+1}}\,D_{2k+1}(k)=-\frac{f(\epsilon)}{\kappa} {\cal D}_{2k},\qquad {\cal D}_{2k}:= \bl(\frac{2}{\psi''(\epsilon)}\br)^k D_{2k+1}(k).
\ee

Hence, at coalescence we obtain the expansion in the following form:
\begin{theorem}$\!\!\!.$\  The expansion of $F_1(\lambda;x)$ as $\lambda\to+\infty$ at coalescence $(\alpha=\epsilon)$ is given by
\bee\label{e28}
F_1(\lambda;x)\sim\frac{G(\lambda)e^{-\lambda\psi(\epsilon)}}{2x}\,f(\epsilon)\bl\{1-\frac{2}{\sqrt{2\pi|\psi''(\epsilon)|}} \sum_{k\geq 0} {\cal D}_{2k}\frac{(\fs)_k}{\lambda^{k+\frac{1}{2}}}\br\}.
\ee 
where the coefficients ${\cal D}_{2k}$ are
\[{\cal D}_0=F_1-\f{1}{6}\Psi_3,\]
\[{\cal D}_2=\frac{1}{\psi''(\epsilon)}\{\f{1}{3}F_3-\f{1}{2}F_2\Psi_3+\f{1}{4}F_1(\f{5}{3}\Psi_3^2-\Psi_4)-\f{1}{4}(\f{35}{54}\Psi_3^3-\f{5}{6}\Psi_3\Psi_4+\f{1}{5}\Psi_5)\},
\]
\[{\cal D}_4=\frac{1}{(\psi''(\epsilon))^2}\{\f{1}{30}F_5-\f{5}{36}F_4\Psi_3+\f{5}{36}F_3(\f{7}{3}\Psi_3^2-\Psi_4)+\f{1}{12}F_2(-\f{35}{6}\Psi_3^3+\f{35}{6}\Psi_3\Psi_4-\Psi_5)
\]
\bee\label{e28a}
\hspace{1.5cm}+\f{1}{36}F_1(\f{385}{24}\Psi_3^4-\f{105}{4}\Psi_3^2\Psi_4+7\Psi_3\Psi_5+\f{35}{8}\Psi_4^2-\Psi_6)
-\f{1}{48}(\f{1001}{108}\Psi_3^5-\f{385}{18}\Psi_3^3\Psi_4\]
\[+\f{35}{4}\Psi_3(\Psi_4^2+\f{4}{5}\Psi_3\Psi_5)-\f{7}{3}(\Psi_4\Psi_5+\f{2}{3}\Psi_3\Psi_6)+\f{4}{21}\Psi_7)\}.
\ee
Here $\Psi_k\equiv \psi^{(k)}(\epsilon)/\psi''(\epsilon)$ and $F_k\equiv f^{(k)}(\epsilon)/f(\epsilon)$, where $f(t)=t^{a-1}(t-1)^{c-a-1}$.
\end{theorem}

\vspace{0.6cm}

\begin{center}
{\bf 4. \  The expansion of $F_m(\lambda;x)$}
\end{center}
\setcounter{section}{4}
\setcounter{equation}{0}
\renewcommand{\theequation}{\arabic{section}.\arabic{equation}}
When the parameter $b=m$, where integer $m\geq 2$, we use the contiguous relation \cite[(15.5.11)]{DLMF}
\[F\left(\!\!\begin{array}{c}\alpha+1, \beta\\ \gamma\end{array}\!;x\!\right)=\frac{(2\alpha-\gamma+(\beta-\alpha)x)}{\alpha(1-x)}\,F\left(\!\!\begin{array}{c}\alpha, \beta\\ \gamma\end{array}\!;x\!\right)+\frac{\gamma-\alpha}{\alpha(1-x)}\,F\left(\!\!\begin{array}{c}\alpha-1, \beta\\ \gamma\end{array}\!;x\!\right).\]
Then, if we define
\[A_m:=-\frac{\lambda}{m(1-x)}\bl\{1-\epsilon x+\frac{c-2m+(m-a)x}{\lambda}\br\},\qquad B_m:=\frac{\lambda}{m(1-x)}\br\{1+\frac{c-m}{\lambda}\br\},\]
we find that
\bee\label{e31}
F_2(\lambda;x)=A_1 F_1(\lambda;x)+B_1.
\ee
From the expansions for $F_1(\lambda;x)$ obtained in Sections 2 and 3 it is then possible to determine the expansion of $F_2(\lambda;x)$. 
Continuation of this process leads to
\[F_3(\lambda;x)=A_2 F_2(\lambda;x)+B_2 F_1(\lambda;x)=(A_1 A_2+B_2) F_1(\lambda;x)+A_2B_1\]
and so on.

An obvious drawback to this approach is the fact that, in the case of (\ref{e31}), we have to evaluate the difference between two terms each multiplied by the large parameter $\lambda$. This problem becomes {\it a fortiori} more acute for higher values of $m$. 

An alternative procedure is to employ an integration by parts. In the case of $F_2(\lambda;x)$ we have
\[F_2(\lambda;x)=\frac{G(\lambda)}{2\pi i x^2} \int_0^{(1+)}\frac{e^{-\lambda\psi(t)}f(t)}{(t-\alpha)^2}\,dt
=\frac{G(\lambda)}{2\pi i x^2} \int_0^{(1+)} \frac{[e^{-\lambda\psi(t)}f(t)]'}{t-\alpha}\,dt\]
since the integrated part vanishes on account of the fact that the integrand near $t=0$ is controlled by $t^{a+\epsilon\lambda-1}$. Then we have the integral (\ref{e21}) with $e^{-\lambda\psi(t)} f(t)$ replaced by
its derivative. Repetition of the argument employed in Section 2 then yields
\[F_2(\lambda;x)\sim -\frac{G(\lambda)}{2x^2}\bl\{[e^{-\lambda\psi(t)}f(t)]'_{t=\alpha} \,\mbox{erfc} (\pm\lambda^{\frac{1}{2}}p)+\frac{e^{-\lambda\psi(\epsilon)}}{\pi}\sum_{k\geq0}B_{2k}\frac{\g(k+\fs)}{\lambda^{k+\frac{1}{2}}}\br\},\]
where
\[B_0=\frac{f'(\epsilon)}{\kappa (\alpha-\epsilon)}\mp \frac{e^{\lambda\psi(\alpha)}}{p}\,[e^{-\lambda\psi(t)}f(t)]'_{t=\alpha}\]
and the choice of signs is as indicated in Theorem 1. The coefficients $B_{2k}$ ($k\geq 1$) are evaluated 
using (\ref{e25w}) for $C_{2k}$ with $f(t)$ replaced by $f'(t)$ and $\lambda\psi'(t)f(t)$. The general case (with $m$ finite) follows from
\[F_m(\lambda;x)=\frac{(-)^m G(\lambda)}{2\pi ix^m \g(m)} \int_0^{(1+)} \frac{[e^{-\lambda\psi(t)}f(t)]^{(m-1)}}{t-\alpha}\,dt.\]

\vspace{0.6cm}

\begin{center}
{\bf 5. \  Numerical results}
\end{center}
\setcounter{section}{5}
\setcounter{equation}{0}
\renewcommand{\theequation}{\arabic{section}.\arabic{equation}}
In Table 1 we display the values\footnote{We have adopted the convention in the tables of writing $x(y)$ for $x\times 10^y$.}  of the coefficients $d_{2k}$ computed from (\ref{e27a}) for $0\leq k\leq 5$ and different values of $x$ when $\epsilon=2$, $a=\fs$ and $c=2$. Coalescence in this case occurs when $x=0.50$. Table 2 displays the values of the absolute relative error in the computation of $F_1(\lambda;x)$ from (\ref{e27b}) when $\epsilon=2$ as a function of the truncation index $k\leq M$ and different values of $\lambda$ and $x$ with the values of $a$ and $c$ specified above.
\begin{table}[th]
\caption{\footnotesize{Values of the coefficients $d_{2k}$ in the case $b=1$ for different $x$ when $\epsilon=2$, $a=\fs$ and $c=2$}} \label{t1}
\begin{center}
\begin{tabular}{|l||l|l|l|l|}
\hline
&&&&\\[-0.3cm]
\mcol{1}{|c||}{$k$} & \mcol{1}{c|}{$x=0.30$} & \mcol{1}{c|}{$x=0.45$} & \mcol{1}{c|}{$x=0.55$}
& \mcol{1}{c|}{$x=0.70$}\\
\hline
&&&&\\[-0.3cm]
0 & $-0.94304503$ & $-1.03364259$ & $-1.08679035$ & $-1.16314077$ \\
1 & $+2.22692591(-1)$ & $+2.22370609(-1)$ & $+2.19071374(-1)$ & $+2.10564977(-1)$ \\
2 & $-1.70235645(-2)$ & $-1.44683556(-2)$ & $-1.32197968(-2)$ & $-1.11672668(-2)$ \\
3 & $-2.01563398(-3)$ & $-2.01674811(-3)$ & $-1.87396028(-3)$ & $-1.50090398(-3)$ \\
4 & $+3.52259554(-4)$ & $+2.80957497(-4)$ & $+2.53010021(-4)$ & $+2.26973127(-4)$ \\
5 & $+4.16738056(-5)$ & $+4.39726773(-5)$ & $+4.00222956(-5)$ & $+2.82024739(-5)$ \\
[.1cm]\hline\end{tabular}
\end{center}
\end{table}
\begin{table}[th]
\caption{\footnotesize{Values of the absolute relative error in the computation of $F_1(\lambda;x)$ from (\ref{e27b}) as a function of the truncation index $M$ and different values of $\lambda$ and $x$ when $\epsilon=2$, $a=\fs$ and $c=2$}} \label{t2}
\begin{center}
\begin{tabular}{|l||l|l|l|l|}
\hline
\mcol{1}{|c}{} & \mcol{4}{c|}{$\lambda=50$}\\
\mcol{1}{|c||}{$M$} & \mcol{1}{c|}{$x=0.30$} & \mcol{1}{c|}{$x=0.45$} & \mcol{1}{c|}{$x=0.55$}
& \mcol{1}{c|}{$x=0.70$}\\
\hline
&&&&\\[-0.3cm]
0 & $2.306(-03)$ & $5.802(-04)$ & $9.799(-05)$ & $2.506(-08)$ \\
1 & $5.331(-06)$ & $1.143(-06)$ & $1.790(-07)$ & $4.201(-11)$ \\
2 & $3.099(-08)$ & $7.831(-09)$ & $1.248(-09)$ & $2,655(-13)$ \\
3 & $3.878(-10)$ & $7.820(-11)$ & $1.207(-11)$ & $2.872(-15)$ \\
4 & $3.989(-12)$ & $1.067(-12)$ & $1.665(-13)$ & $3.104(-17)$ \\
5 & $9.809(-14)$ & $1.916(-14)$ & $2.964(-15)$ & $7.326(-19)$ \\
[.1cm]\hline
\mcol{1}{|c}{} & \mcol{4}{c|}{$\lambda=100$}\\
\mcol{1}{|c||}{$M$} & \mcol{1}{c|}{$x=0.30$} & \mcol{1}{c|}{$x=0.45$} & \mcol{1}{c|}{$x=0.55$}
& \mcol{1}{c|}{$x=0.70$}\\
\hline
&&&&\\[-0.3cm]
0 & $1.103(-03)$ & $2.426(-04)$ & $1.861(-05)$ & $1.451(-12)$ \\
1 & $1.270(-06)$ & $2.378(-07)$ & $1.692(-08)$ & $1.212(-15)$ \\
2 & $3.724(-09)$ & $8.218(-10)$ & $5.947(-11)$ & $3.862(-18)$ \\
3 & $2.304(-11)$ & $4.055(-12)$ & $2.843(-13)$ & $2.066(-20)$ \\
4 & $1.206(-13)$ & $2.811(-14)$ & $1.993(-15)$ & $1.136(-22)$ \\
5 & $1.455(-15)$ & $2.477(-16)$ & $1.741(-17)$ & $1.316(-24)$ \\
[.1cm]\hline
\end{tabular}
\end{center}
\end{table}
\begin{table}[h]
\caption{\footnotesize{Values of the coefficients ${\cal D}_{2k}$ at coalescence in the case $b=1$ for $x=0.50$ when $\epsilon=2$, $a=\fs$ and $c=2$. The corresponding values of $d_{2k}$ are also shown.}} 
\begin{center}
\begin{tabular}{|l||l|l|}
\hline
&&\\[-0.3cm]
\mcol{1}{|c||}{$k$} & \mcol{1}{c|}{${\cal D}_{2k}$} & \mcol{1}{c|}{$d_{2k}$} \\
\hline
&&\\[-0.3cm]
0 & $+0.75000000$ & $-1.06066017$ \\
1 & $-1.56250000(-1)$ & $+2.20970869(-1)$ \\
2 & $+9.76562500(-2)$ & $-1.38106793(-2)$ \\
[.1cm]\hline\end{tabular}
\end{center}
\end{table}
\begin{table}[h]
\caption{\footnotesize{Values of the absolute relative error in the computation of $F_1(\lambda;x)$ at coalescence ($x=0.50$, $\epsilon=2$) from (\ref{e28}) as a function of the truncation index $M$ and different values of $\lambda$ when $a=\fs$ and $c=2$}} 
\begin{center}
\begin{tabular}{|l||l|l|l|}
\hline
&&&\\[-0.3cm]
\mcol{1}{|c||}{$M$} & \mcol{1}{c|}{$\lambda=50$} & \mcol{1}{c|}{$\lambda=100$} & \mcol{1}{c|}{$\lambda=150$}\\
\hline
&&&\\[-0.3cm]
0 & $2.826(-04)$ & $9.620(-05)$ & $5.151(-05)$ \\
1 & $5.346(-07)$ & $9.050(-08)$ & $3.229(-08)$ \\
2 & $3.724(-09)$ & $3.182(-10)$ & $7.582(-11)$ \\
[.1cm]\hline\end{tabular}
\end{center}
\end{table}

In the case of coalescence (when $\epsilon=2$) we present the coefficients ${\cal D}_{2k}$ defined in (\ref{e28a}) for $0\leq k\leq 2$ in Table 3. From (\ref{e28b}) we have $d_{2k}=-2^{1/2} {\cal D}_{2k}$ in this case. It is seen that the values of $d_{2k}$ `sit' in between the values corresponding to $x=0.45$ and $x=0.55$ in Table 1. Finally, Table 4 shows an example of the absolute relative error in the computation of $F_1(\lambda;\fs)$ at coalescence when $\epsilon=2$ for different $\lambda$ and truncation index $M$.
\vspace{0.6cm}

\begin{center}
{\bf 6. \  Concluding remarks}
\end{center}
\setcounter{section}{6}
\setcounter{equation}{0}
\renewcommand{\theequation}{\arabic{section}.\arabic{equation}}
We have derived a uniform expansion for the hypergeometric function $F_1(\lambda;x)=F(a+\epsilon\lambda,1;c+\lambda;x)$ for $\lambda\to+\infty$ when $\epsilon>1$ and $0<x<1$. The resulting expansion is presented in Theorem 1 and the corresponding expansion at coalescence when $\epsilon x=1$ is stated in Theorem 2. 

The expansion of the function when the second numerator parameter equals $m$, with integer $m\geq 2$, is covered by a recurrence relation in Section 4. Clearly this is not an ideal situation. It is hoped that the expansion for $m\geq 2$ can be elaborated from an integral representation of the type (\ref{e21}). In addition, the case of non-integer $m$ also needs to be considered as outlined in \cite{B}.

Another limitation inherent in the results is the expansion of $F_1(\lambda;x)$ in the neighbourhood of $\epsilon x=1$.
Although it is possible to use the expansion in Theorem 1 as $\epsilon x\to 1$, this would necessitate using increasing precision to deal with the removable singularity present in the coefficients $d_{2k}$ at coalescence.

\vspace{0.6cm}

\begin{center}
{\bf Appendix: An alternative form of expansion }
\end{center}
\setcounter{section}{1}
\setcounter{equation}{0}
\renewcommand{\theequation}{\Alph{section}.\arabic{equation}}
An alternative form of the expansion (\ref{e27b}) can be obtained by noting that the contribution
to $F_1(\lambda;x)$ resulting from the sum involving the coefficients $b_{2k}$ can be written as
\[
-\frac{G(\lambda) e^{-\lambda\psi(\alpha)}}{2x}\,f(\alpha)S(\lambda p^2),\qquad S(\lambda p^2):=\frac{e^{-\lambda p^2}}{\sqrt{\pi}}  \sum_{k\geq 0} (-)^k(\fs)_k (\lambda p^2)^{-k-\frac{1}{2}}.
\]
From \cite[(7.12.1)]{DLMF} it is seen that $S(\lambda p^2)$ is the asymptotic expansion of $\mbox{erfc}\,(\lambda^\frac{1}{2}p)$ as $\lambda p^2\to+\infty$.

If we now define
\[\mbox{Erfc}\,(\pm x):=\mbox{erfc}\,(\pm x)\mp\frac{e^{-x^2}}{\sqrt{\pi} x} \sum_{k\geq 0} (\fs)_k (-x^2)^{-k},\]
with the usual interpretation of an asymptotic series, we obtain the alternative forms of the expansion given by
\bee\label{a1}
F_1(\lambda;x)\sim G(\lambda)\bl\{\frac{e^{-\lambda\psi(\alpha)} f(\alpha)}{2x}\, \mbox{Erfc}\,(\pm\lambda^\fr p)+\frac{e^{-\lambda\psi(t_s)}}{\sqrt{2\pi|\psi''(\epsilon)|}}\,\frac{f(\epsilon)}{(1-\epsilon x)}\sum_{k\geq 0} C_{2k} \frac{(\fs)_k}{\lambda^{k+1/2}}\br\}.
\ee

\end{document}